\documentclass[12pt, notitlepage]{amsart}
\usepackage{latexsym, amsfonts, amsmath, amssymb, amsthm, cite}

\newtheorem{lemma}{Lemma}

\newtheorem{theorem}{Theorem}

\newcommand{\assign}{:=}

\usepackage{graphicx}
\usepackage{mathrsfs}

\begin{document}

\title{Continuous lower bounds for moments of  zeta and $L$-functions}
\author{Maksym Radziwi\l\l}
\address{Department of Mathematics \\ Stanford University \\
450 Serra Mall, Bldg. 380\\ Stanford, CA 94305-2125}
\email{maksym@stanford.edu}
\author{Kannan Soundararajan}
\address{Department of Mathematics \\ Stanford University \\
450 Serra Mall, Bldg. 380\\ Stanford, CA 94305-2125}
\email{ksound@math.stanford.edu}
\thanks{The first author is partially supported by a NSERC PGS-D award.  The 
second author is partially supported by NSF grant DMS-1001068}

\dedicatory{In memoriam Professor K. Ramachandra (1933-2011)}

\subjclass[2010]{Primary: 11M06, Secondary: 11M50}

\begin{abstract}
  We obtain lower bounds of the correct order of magnitude for the $2 k$-th
  moment of the Riemann zeta function for all $k \ge 1$. Previously such lower bounds 
  were known only for rational values of $k$, with the bounds depending on the 
  height of the rational number $k$.  Our new bounds are continuous in $k$, and thus 
  extend also to the case when $k$ is irrational.    The method is a refinement 
  of an approach of Rudnick and Soundararajan, and applies also to moments of $L$-functions 
  in families.
  \end{abstract}

\maketitle

\section{Introduction}

\noindent Ramachandra \cite{Ram1, Ram2} established that the moments 
\[ 
M_k (T) = \int_0^{T} | \zeta ( \tfrac{1}{2} + i t) |^{2 k} dt
\]
satisfy the lower bound 
$$
M_k(T) \ge (c_k +o(1)) T (\log T)^{k^2}
$$ 
for natural numbers $2k$.   Here $c_k$ is a positive constant, and the estimate holds as $T$ tends to infinity.  Ramachandra's result was extended by Heath-Brown \cite{HeBr1} to the case when $k$ is any positive rational number.   The constant $c_k$ in Heath-Brown's result depends on the height of 
the rational number $k$, and does not vary continuously 
with $k$.  Thus when $k$ is irrational, only the weaker result, due to Ramachandra \cite{Ram3},  
$M_k(T) \gg T (\log T)^{k^2} (\log \log T)^{-k^2}$ is known.  If the truth of the Riemann Hypothesis 
is assumed, then Ramachandra \cite{Ram2, Ram4} and independently Heath-Brown \cite{HeBr1} showed that it is possible to get the stronger bound $M_k(T) 
\gg_k T(\log T)^{k^2}$ for all positive real $k$.  

The lower bounds discussed above are essentially of the right order of magnitude for $M_k(T)$.  
A folklore conjecture states that $M_k(T) \sim C_k T (\log T)^{k^2}$ for a positive 
constant $C_k$ and all positive real numbers $k$.   This is known for $k=1$ and $2$ (see 
Chapter VII of \cite{Tit}), but remains open for all other values.  A more precise conjecture, 
predicting the value for $C_k$, has been proposed by Keating and Snaith \cite{KS1}.  When 
$0\le k\le 2$, Heath-Brown \cite{HeBr1} showed (assuming RH) that $M_k(T) \ll_k T(\log T)^{k^2}$, 
and this has recently been extended by Radziwi\l\l \cite{Rad} (also on RH) to cover the range $k \le 2.18$.  
For all positive $k$, Soundararajan \cite{Sou1} has shown (on RH) that $M_k(T) \ll_{k,\epsilon} T(\log T)^{k^2+\epsilon}$.

 We now return to the problem of lower bounds for moments which is the focus of this 
 paper.  Recently, Rudnick and Soundararajan \cite{RuS1, RuS2} 
 developed a method for establishing lower bounds of the conjectured order of magnitude 
 for rational moments of $L$-functions varying in families.   As with Heath-Brown's result 
 for $\zeta(s)$, these lower bounds do not vary continuously with the parameter $k$, 
 and thus one fails to get lower bounds of the right order of magnitude for irrational moments.  
  In this paper, we show how the Rudnick-Soundararajan approach may be extended 
  so as to obtain lower bounds that vary continuously with the moment parameter $k$.  
  In particular, we thus obtain lower bounds of the right order of magnitude for all real $k\ge 1$.  
  Since this is new already for $\zeta(s)$, we present the proof for that case, and sketch (in section 5)  
  the modifications for $L$-functions.

\begin{theorem}   For any real number $k >1$, and all large $T$ we have 
$$ 
M_k(T) \ge e^{-30 k^4} T (\log T)^{k^2}. 
$$ 
\end{theorem}

We have made no effort to obtain the best possible constant in our Theorem, and 
with more work it may be possible to obtain a substantially better value.   When 
$k$ is a natural number, Conrey and Ghosh \cite{CG} gave the elegant lower bound 
$$ 
M_k(T) \ge (1+o(1)) T \sum_{n\le T} \frac{d_k(n)^2}{n} 
\sim T \frac{(\log T)^{k^2}}{\Gamma(k^2+1)} \prod_p \Big(1-\frac 1p\Big)^{k^2} 
\Big(1+\sum_{a=1}^{\infty} \frac{d_k(p^a)^2}{p^a}\Big), 
$$ 
and the best known lower bound (see \cite{Sou2}) is twice as large.   These 
bounds are of the rough shape  $M_k(T) \gg k^{-2k^2 (1+o_k(1))} T (\log T)^{k^2}$ for 
large natural numbers $k$.  The  constant in the conjectured asymptotic formula for $M_k(T)$ is roughly of size $k^{-k^2 (1+o(1))}$ (see, for example, the discussion in \cite{CoGo}).

As in \cite{RuS1}, the method presented here applies to any family
of $L$-functions for which slightly more than the first moment can be
understood, and produces good lower bounds for all moments larger than the first.
There still remains the question of finding lower bounds when $0< k <1$.   Recently 
Chandee and Li \cite{CL} developed a method for producing good lower bounds for  moments of 
$L$-functions when $0<k<1$ is rational.   In forthcoming work, Radziwi\l\l describes an alternative approach 
which establishes that $M_k(T) \gg_k T(\log T)^{k^2}$ for all real $0<k <1$, but his approach is 
special to the $t$-aspect and does not adapt to central values of $L$-functions.

\section{Proof of Theorem 1}

\noindent For any real number $0< \alpha\le 1$ define the {\sl Sylvester sequence} $s_n=s_n(\alpha)$ 
as follows:  $s_1$ is the least integer strictly larger than $1/\alpha$, and 
given $s_1$, $\ldots$, $s_n$ take $s_{n+1}$ to be the least integer 
strictly larger than $(\alpha -\frac 1{s_1}-\ldots -\frac1{s_n})^{-1}$.   
Thus $\alpha = \sum_{n=1}^{\infty} \frac{1}{s_n(\alpha)}$.   When $\alpha =1$ 
we obtain the sequence $2$, $3$, $7$, $43$, $\ldots$, where the 
next entry in the sequence is obtained by multiplying the all preceding entries and 
adding $1$, or equivalently by the recurrence $s_{n+1}(1) = s_n(1) (s_n(1)-1) +1$.   
For any $0 < \alpha \le 1$, we may easily check that $s_1(\alpha) \ge 2$ and 
$s_{n+1}(\alpha) \ge s_n(\alpha)(s_n(\alpha)-1) + 1$, so that $s_n(\alpha) \ge s_n(1)$ 
for all $n$.   The Sylvester sequence grows very rapidly: from the recurrence for $s_{n+1}(1)$ 
we see easily that $s_n(1) \ge (n-1)! + 1$, and indeed one can show that $s_n(1)$ grows 
doubly exponentially.

For any $k>1$, we denote by $a_\ell$ the Sylvester sequence for $1-\frac 1k$, 
and by $b_\ell$ the Sylvester sequence for $1$.  Thus $a_\ell \ge b_\ell \ge (\ell-1)!+1$, 
and 
$$ 
\sum_{\ell=1}^{\infty} \frac1{a_\ell} = 1-\frac 1k, \qquad\text{and}\qquad 
\sum_{\ell=1}^{\infty} \frac1{b_\ell} =1. 
$$ 
Let  $T$ be large, and set $T_0= T^{1-\vartheta}$, for a small positive parameter $\vartheta$; for example, we may simply
take $\theta = 1/100$ below.  We define the Dirichlet polynomials
$$ 
A_\ell (s) = \sum_{n\le T_0^{1/a_\ell}} \frac{d_{k/a_\ell}(n)}{n^s}, 
\qquad \text{and} \qquad 
B_\ell(s) = \sum_{n\le T_0^{1/b_\ell}} \frac{d_{k/b_\ell}(n)}{n^s}. 
$$
Let $K$ denote a smooth non-negative function compactly supported inside the interval $[\vartheta,1-\vartheta]$,  
with $K(x) \le 1$ for all $x$, and $K(x)=1$ for $x\in [2\vartheta, 1-2\vartheta]$.   For such 
$K$ we have that the Fourier transform ${\hat K}(\xi) = \int_{-\infty}^{\infty} K(x) e^{-ix\xi} dx$ 
satisfies for any non-negative integer $\nu$ 
\begin{equation}
\label{eqn2.0} 
|{\hat K}(\xi) | \ll (1+|\xi|)^{-\nu}. 
\end{equation} 

 Our refinement of the method of Rudnick and Soundararajan is based upon a consideration 
 of the quantity
 \begin{equation} 
 \label{eqn:2.1}
 \mathcal{I}(T) \assign  \int_{0}^{T} K\Big(\frac{t}{T}\Big)  \zeta ( \tfrac{1}{2} +
   i t)  \prod_{\ell=1}^{\infty} A_{\ell} ( \tfrac{1}{2} +i
   t) B_{\ell} ( \tfrac{1}{2} -  i t)  dt. 
  \end{equation} 
  Note that if $\ell$ is sufficiently large, then $A_\ell(s)$ and $B_\ell(s)$ will be identically 
  $1$, and so in the infinite product above only finitely many terms matter.  
 On the one hand, we shall establish a good lower bound for ${\mathcal I}(T)$.  
 
 \begin{lemma} 
 \label{lemma1} 
 With the above notations, we have 
 $$ 
 {\mathcal I}(T) \ge e^{-15 k^3} T (\log T)^{k^2}. 
 $$ 
 \end{lemma} 
 
  On the other hand,  by H{\" o}lder's inequality we have that   
  bounded by 
  \begin{eqnarray*}
{\mathcal I}(T) &\le& \Big(\int_0^{T} K\Big(\frac tT\Big) |\zeta(\tfrac 12+it)|^{2k} dt\Big)^{\frac 1{2k}}  \\
 &&\times  \prod_{\ell =1}^{\infty}
  \Big( \int_0^{T} K\Big(\frac{ t}{T}\Big) |A_\ell (\tfrac 12+it)|^{2a_\ell} dt \Big)^{\frac{1}{2a_\ell}} 
  \Big(\int_0^{T} K\Big(\frac{t}{T}\Big) |B_\ell(\tfrac 12+it)|^{2b_\ell} dt \Big)^{\frac 1{2b_\ell} },  
  \end{eqnarray*}
  and we may work out upper bounds for the terms in the product above. 
  
    \begin{lemma}
 \label{lemma2}  For any integer $a\ge 1$ we have 
 $$
  \int_0^{T} K\Big( \frac tT\Big) 
   \Big| \sum_{n \le T_0^{1/a}} \frac{d_{k/ a} (n)}{n^{\frac{1}{2} + i t}} \Big|^{2 a} dt 
    \le   T \sum_{n\le T_0} \frac{d_k(n)^2}{n} + O(1) \le T (\log T)^{k^2}.
$$
\end{lemma}

From the above two lemmas and our upper bound on ${\mathcal I}(T)$ we obtain that
$$ 
M_k(T) \ge \frac{{\mathcal I}(T)^{2k}}{(T (\log T)^{k^2})^{(2k-1)}} \ge e^{-30k^4} T (\log T)^{k^2}, 
$$ 
proving our Theorem.  

\section{Proof of Lemma 2} 

\noindent  Write 
  \begin{eqnarray*}
    \Big( \sum_{n \le T_0^{1 / a}} \frac{d_{k / a} (n)}{n^s} \Big)^a
    & = & \sum_{n \le T_0} \frac{a(n)}{n^s},
  \end{eqnarray*}
say,   where 
  \[ 
  0 \le a(n) = \sum_{\substack {n = n_1 \ldots n_a\\
       n_1, \ldots, n_r \le T_0^{1 /a}
     }} d_{k / a} (n_1) \cdot \ldots \cdot d_{k / a} (n_a) \text{ }
     \le \text{ } d_k (n) . \]
  Then
  \begin{eqnarray*}
    \int_{- \infty}^{\infty} K\Big(\frac tT\Big) \Big| \sum_{n \le T_0^{1 / a}}
    \frac{d_{k / a} (n)}{n^{\frac{1}{2} + i t}} \Big|^{2a} dt
    & = &  \int_{- \infty}^{\infty} K\Big( \frac tT\Big) \Big| \sum_{n
    \le T_0}\frac{a (n)}{n^{\frac{1}{2} + i t}}
    \Big|^2 dt 
    \\
    &=& T\sum_{n, m \le T_0} \frac{a(m) a(n)}{\sqrt{mn}} 
     \hat{K} \Big( T \log \frac{n}{m} \Big).\\
     \end{eqnarray*}
     The diagonal terms $m=n$ above contribute 
     $$ 
     T{\hat K}(0) \sum_{n\le T_0} \frac {a(n)^2}{n} \le T\sum_{n\le T_0} \frac{d_k(n)^2}{n}.  
     $$ 
         To handle the off-diagonal terms $m\neq n$, note that if $m \neq n \le T_0$ then 
     $|T\log (n/m)| \gg T/T_0 = T^{\vartheta}$ so that by \eqref{eqn2.0} we have 
     $|{\hat K}(T\log (n/m))| \ll T^{-2}$.  Therefore the off-diagonal terms contribute 
    $$ 
    \ll T^{-1} \sum_{m, n \le T_0} \frac{d_k(m)d_k(n)}{\sqrt{mn}} \ll 1.  
    $$ 
    Adding these contributions, we obtain that 
    $$ 
    {\mathcal I}(T) \le T\sum_{n\le T_0} \frac{d_k(n)^2}{n} + O(1).
    $$ 
    Since 
    $$ 
    \sum_{n\le T_0} \frac{d_k(n)^2}{n } 
    \sim \frac{(\log T_0)^{k^2}}{\Gamma(k^2+1)} \prod_p \Big(1-\frac 1p\Big)^{k^2} 
    \Big( 1+ \sum_{a=1}^{\infty} \frac{d_k(p^a)^2}{p^a} \Big) \le (\log T)^{k^2},
    $$ 
    the Lemma follows.

\section{Proof of Lemma 1}

\noindent  For $s= \tfrac 12+ it$ with $\vartheta T \le  t \le  T$ we have
\[
 \zeta (s) = \sum_{n \le T} \frac{1}{n^s} + O(T^{-\frac 12}) . 
\]
Note that, by H{\" o}lder's inequality and Lemma 2,  
\begin{eqnarray*}
\int_0^T K\Big(\frac tT\Big)\!\!\!\!\!\!&&\!\!\!\! \prod_{\ell=1}^{\infty} |A_\ell(\tfrac 12+it)B_\ell(\tfrac 12+it)| dt \\
&\le& T^{\frac{1}{2k}} \prod_{\ell =1}^{\infty} \Big(\int_0^T K \Big( \frac tT\Big) |A_\ell(\tfrac 12+it)|^{2a_\ell} 
dt \Big)^{\frac 1{2a_\ell}} \Big(\int_0^T K\Big(\frac tT\Big) |B_\ell(\tfrac 12+it)|^{2b_\ell}\Big)^{\frac{1}{2b_\ell}} \\
&\ll& T^{1+\epsilon}. 
\end{eqnarray*} 
Therefore 
\begin{equation}
\label{4.1}
{\mathcal I}(T) = \int_0^{T}K\Big(\frac tT\Big) \sum_{n \le T } \frac{1}{n^{\frac{1}{2} 
+   i t}}  \prod_{\ell=1}^{\infty} A_{\ell} ( \tfrac{1}{2} + i
   t) B_{\ell} ( \tfrac{1}{2} - i t) dt + O (T^{\frac 12+\epsilon}). 
  \end{equation} 
  
  Write 
  $$ 
  \sum_{n\le T} \frac 1{n^s} \prod_{\ell=1}^{\infty} A_\ell(s) =\sum_{n} \frac{\alpha(n)}{n^s}, 
  \qquad \text{and} \qquad \prod_{\ell=1}^{\infty} B_\ell(s) = \sum_{n} \frac{\beta(n)}{n^s}. 
  $$ 
  Then both $\alpha(n)$ and $\beta(n)$ are non-negative and bounded above by $d_k(n)$.  
  Moreover $\alpha(n) =0$ if $n> T^2  (> T T_0^{\sum_\ell 1/a_\ell})$ and $\beta(n)=0$ if 
  $n> T_0 =T_0^{\sum_\ell 1/b_\ell}$.   From \eqref{4.1} we obtain that 
  $$ 
  {\mathcal I}(T) = \sum_{m, n}\frac{\alpha(m)\beta(n)}{\sqrt{mn}} T{\hat K}\Big(T\log \frac nm\Big) 
  + O(T^{\frac 12+\epsilon}). 
  $$ 
  If $m\neq n$ and $n\le T_0$ then $|T\log (n/m)| \gg T/T_0 = T^{\vartheta}$ 
  and so, using \eqref{eqn2.0}, the off-diagonal terms above contribute 
$$ 
\ll T^{-2} \sum_{m\le T^2} \sum_{n\le T_0} \frac{d_k(m)d_k(n)}{\sqrt{mn}} \ll 1. 
$$ 
Thus 
\begin{equation} 
\label{4.2} 
{\mathcal I}(T) = T {\hat K}(0) \sum_{n} \frac{\alpha(n)\beta(n)}{n} + O(T^{\frac 12+\epsilon}).
\end{equation}

 Let $A_0$, $A_\ell$, $B_\ell$ ($\ell \ge 1$) denote parameters all larger than $1$, 
 to be chosen later.  Set $\alpha_0 = A_0/\log T_0$, and for $\ell \ge 1$, $\alpha_\ell = A_\ell/\log T_0$ and $\beta_\ell = B_\ell/\log T_0$.  From the definition, we see that $\alpha(n)\beta(n)$ equals the  
 sum of the 
 quantity $\prod_{\ell \ge 1} d_{k/a_\ell}(m_\ell) d_{k/b_\ell}(n_\ell)$ over all possible ways of writing $n= m_0 \prod_{\ell \ge 1} m_\ell= \prod_{\ell \ge 1} n_\ell$ 
 with $m_0 \le T$, $m_\ell \le T_0^{1/a_\ell}$ and $n_\ell \le T_0^{1/b_\ell}$.  
 Next note that 
 $$ 
 m_0^{-\alpha_0} \prod_{\ell \ge 1} m_\ell^{-\alpha_\ell} n_\ell^{-\beta_\ell } 
 - e^{-A_0} - \sum_{\ell \ge 1} \Big( e^{-A_\ell/a_\ell} + e^{-B_\ell/b_\ell}\Big) 
 $$ 
 is always less than $1$, and is less than $0$ if $m_0 >T (>T_0)$ or if any 
 $m_\ell > T_0^{1/a_\ell}$ or if any $n_\ell > T_0^{1/b_\ell}$.   Therefore we see that 
$\alpha(n) \beta(n)$ is at least as large as 
 \begin{align*}
\sum_{n=m_0\prod_\ell m_\ell =\prod_\ell n_\ell} 
& \prod_\ell d_{k/a_\ell}(m_\ell) d_{k/b_\ell}(n_\ell) 
 \Big(  m_0^{-\alpha_0} \prod_{\ell \ge 1} m_\ell^{-\alpha_\ell} n_\ell^{-\beta_\ell} 
 - e^{-A_0} - \sum_{\ell \ge 1} \Big( e^{-A_\ell/a_\ell} + e^{-B_\ell/b_\ell}\Big) \Big)\\
 &=f(n) - \Big( e^{-A_0} + \sum_{\ell \ge 1} \Big( e^{-A_\ell/a_\ell} + e^{-B_\ell/b_\ell}\Big) \Big) d_k(n)^2, 
 \end{align*}
 say.  Above, we see that $f(n)$ is a multiplicative function of $n$, and for a prime number $p$ we have 
 (adopting the convention $a_0=k$) 
 \begin{equation}
 \label{4.3}
 f(p) = \sum_{j=0}^{\infty} \sum_{\ell =1}^{\infty} \frac{k^2}{a_j b_\ell} p^{-\alpha_j - \beta_\ell}. 
 \end{equation}

Restricting attention to square-free numbers $n$ that are composed only 
of prime factors below $T_0$, we see from the above remarks that 
 \begin{align*}
 \sum_{n} \frac{\alpha(n) \beta(n)}{n}& \ge \sum_{{p|n \implies p\le T_0}} 
 \frac{\mu(n)^2}{n} \Big(f(n) - \Big( e^{-A_0} + \sum_{\ell \ge 1} \Big( e^{-A_\ell/a_\ell} + e^{-B_\ell/b_\ell}\Big) \Big) d_k(n)^2\Big) \\
& = \prod_{p\le T_0} \Big( 1+ \frac{f(p)}{p} \Big)  - \Big( e^{-A_0} + \sum_{\ell \ge 1} \Big( e^{-A_\ell/a_\ell} + e^{-B_\ell/b_\ell}\Big) \Big) \prod_{p\le T_0} \Big( 1+ \frac{k^2}{p}\Big). 
 \end{align*}
From \eqref{4.3} we see that $f(p) \le k^2$, and since $(1+y)/(1+x) \ge \exp(y-x)$ whenever $0\le y\le x$, we have $(1+f(p)/p) \ge (1+k^2/p) \exp((f(p)-k^2)/p)$.  Thus we obtain that the above is 
\begin{equation} 
\label{4.4}
\ge \prod_{p\le T_0} \Big( 1+\frac{k^2}{p} \Big) 
\Big( \exp\Big(-\sum_{p\le T_0} \sum_{j=0}^{\infty} 
\sum_{\ell=1}^{\infty} 
\frac{k^2}{a_j b_\ell} \Big(\frac 1p -\frac{1}{p^{1+\alpha_j+\beta_\ell}}\Big)\Big) - 
 \Big( e^{-A_0} + \sum_{\ell \ge 1} \Big( e^{-A_\ell/a_\ell} + e^{-B_\ell/b_\ell}\Big) \Big) \Big).
 \end{equation}
 
 If $x$ is large and $\alpha > (\log x)^{-\frac 12}$ then 
 \begin{align*}
 \sum_{p\le x} \Big(\frac 1p - \frac{1}{p^{1+\alpha}} \Big) 
& \le \sum_{p\le x}\Big(  \log \Big(1-\frac 1p\Big)^{-1} - \log \Big( 1-\frac{1}{p^{1+\alpha}}\Big) \Big)
\\
& = \log \log x + \gamma - \log \zeta(1+\alpha) +o(1) \le \log (\alpha \log x) + \gamma +o(1). 
\end{align*} 
On the other hand, if $(\log x)^{-1} < \alpha \le (\log x)^{-\frac 12}$ then 
\begin{align*} 
 \sum_{p\le x} \Big(\frac 1p - \frac{1}{p^{1+\alpha}} \Big) 
& \le \sum_{p\le e^{1/\alpha}} \frac{\alpha \log p}{ p}  + \sum_{e^{1/\alpha} \le p \le x} \frac 1p  \le 1 + \log (\alpha \log x) + o(1). 
\end{align*}
The second bound works for all $\alpha \ge (\log x)^{-1}$, and using it in \eqref{4.4} and 
substituting that back into \eqref{4.2}, we conclude that 
\begin{align*} 
 {\mathcal I}(T)& \ge (1+o(1)) T {\hat K}(0) \prod_{p\le T_0} \Big( 1+ \frac{k^2}{p}\Big) 
 \\ 
 &\times \Big( 
\exp \Big( - \sum_{j=0}^{\infty} \sum_{\ell = 1}^{\infty} 
\frac{k^2}{a_j b_\ell} (\log (A_j+B_\ell) +1)\Big) -e^{-A_0} - \sum_{\ell=1}^{\infty} (e^{-A_\ell/a_\ell} +e^{-B_\ell/b_\ell})\Big). 
\end{align*} 
 
  We now choose $A_0 = 20 k^3$, 
and for $\ell \ge 1$ choose $A_\ell = 20 k^3 a_\ell^2$ and $B_\ell = 20 k^3 b_\ell^2$.  Then 
we see easily that 
$$ 
e^{-A_0} + \sum_{\ell=1}^{\infty} (e^{-A_\ell/a_\ell} + e^{-B_\ell/b_\ell}) 
\le 2 e^{-20 k^3}.
$$  
Further, using that $\log (a_j^2+b_\ell^2) \le \log (1+a_j^2) + \log (1+b_\ell^2)$, 
\begin{align*} 
\sum_{j=0}^{\infty} \sum_{\ell=1}^{\infty} \frac{k^2}{a_j b_\ell} (\log (A_j+B_\ell) +1)
&= k^2 + \sum_{\ell=1}^{\infty} \frac{k}{b_\ell} \log (20k^3 (1+b_\ell^2)) + 
\sum_{j=1}^{\infty} \sum_{\ell =1}^{\infty} 
\frac{k^2}{a_j b_\ell} \log (20k^3(a_j^2+b_\ell^2)) \\ 
&\le k^2 \Big(1 + \log (20k^3) + \sum_{\ell = 1}^{\infty} \frac{\log (1+b_\ell^2)}{b_\ell}  
+ \sum_{j=1}^{\infty} \frac{\log (1+a_j^2)}{a_j}\Big). 
\end{align*}
Now by a calculation we find that $\sum_{\ell=1}^{\infty} \log (1+b_\ell^2)/b_\ell < 5/2$, 
and since $a_j \ge b_j$ for $j\ge 1$ and $\log (1+x^2)/x$ is decreasing for $x\ge 2$, 
it also follows that $\sum_{j=1}^{\infty} \log (1+a_j^2)/a_j < 5/2$.   Thus we conclude that 
$$ 
\sum_{j=0}^{\infty} \sum_{\ell=1}^{\infty} \frac{k^2}{a_j b_\ell} (\log (A_j+B_\ell) +1) 
\le k^2 (6 + \log (20k^3)) \le 10k^3.
$$ 
From these estimates, we obtain that for any $k> 1$ 
$$ 
{\mathcal I}(T) \ge (1+o(1)) T {\hat K}(0) e^{-12 k^3} \prod_{p\le T_0} \Big( 1+\frac{k^2}{p}\Big). 
$$ 

Restricting attention to $p>k$, and since $(1+x) \ge \exp(x-x^2/2)$ for $0\le x\le 1$, we 
have 
$$ 
\prod_{p\le T_0} \Big(1+ \frac{k^2}{p} \Big) 
\ge \exp\Big( \sum_{k<p \le T_0} \Big(\frac{k^2}{p} - \frac{k^4}{2p^2} \Big) \Big) 
\ge \exp\Big( -2 k^3 + \sum_{p\le T_0} \frac{k^2}{p} \Big) \ge e^{-2k^3} (\log T)^{k^2},
$$ 
since for large $T$, $\sum_{p\le T_0} 1/p  = \log \log T_0 + B_1 + o(1) \ge \log \log T$, 
where $B_1= 0.2614 \ldots$.   Using this in our bound for ${\mathcal I}(T)$, and since
$\hat{K}(0) \geq 1 - 4\theta \geq 3/5$ if $\theta < 1/10$, the Lemma follows at once.

\section{Remarks} 

\noindent We may explain the success of our method as follows.  Write 
$$ 
M_k(T) = \int_0^T \zeta(\tfrac 12+it) \prod_{\ell=1}^{\infty} \zeta(\tfrac 12+it)^{k/a_\ell} 
\zeta(\tfrac 12-it)^{k/b_\ell} dt. 
$$ 
Our idea is then to ``approximate" $\zeta(\tfrac 12+it)^{k/a_\ell}$ by $A_\ell (\frac 12+it)$ 
and $\zeta(\tfrac 12-it)^{k/b_\ell}$ by $B_\ell (\tfrac 12-it)$.  Note that 
$A_\ell$ and $B_\ell$ are short Dirichlet polynomials, with diminishing length as $\ell$ increases. 
This permits the evaluation of the quantity ${\mathcal I}(T)$.  On the other hand, we 
would expect that as $\ell$ increases, the terms 
$\zeta(\tfrac 12+it)^{k/a_\ell} \zeta(\tfrac 12-it)^{k/b_\ell}$ make progressively smaller 
impacts on the moment $M_k(T)$, so that approximating these quantities by shorter 
Dirichlet polynomials does not entail too great a loss.   While the Sylvester sequences 
seem a natural choice in this construction, all that we require is the 
convergence of $\sum_{\ell} \log (1+a_\ell)/a_\ell$ and $\sum_{\ell }\log (1+b_\ell)/b_\ell$.

We may easily modify this method to the case of $L$-functions in families.  
For example, if $q$ is a large prime we may start with 
$$
{\mathcal I}(q) = \sum_{\chi \pmod q}^* L(\tfrac 12, \chi) \prod_{\ell =1}^{\infty} A_\ell(\chi) B_\ell (\overline{\chi}), 
$$ 
where the sum is over primitive characters $\chi$, and $A_\ell (\chi) = \sum_{n\le q^{\vartheta/a_\ell}} d_{k/a_\ell}(n)\chi(n)/\sqrt{n}$ 
and $B_\ell(\overline{\chi})= \sum_{n\le q^{\vartheta/b_\ell}} d_{k/b_\ell}(n) \overline{\chi}(n)/\sqrt{n}$.  If $\vartheta$ is small, then $\prod_{\ell} A_\ell(\chi) = \sum_{n} \alpha(n)\chi(n)/\sqrt{n}$ 
and $\prod_{\ell}B_\ell (\overline{\chi}) = \sum_{n} \beta(n) \overline{\chi}(n)/\sqrt{n}$ 
are short sums with $\alpha(n)=\beta(n) =0$ if $n\ge q^{\vartheta}$.   Therefore ${\mathcal I}(q)$ 
behaves like a moment of $L(\tfrac 12, \chi)$ slightly larger than the first, and arguing as 
in \cite{RuS1} we may evaluate ${\mathcal I}(q)$.   Then by applying H{\" o}lder's inequality 
we obtain lower bounds for $\sum_{\chi \pmod q}^* |L(\tfrac 12,\chi)|^{2k}$.   

To take another example, consider 
$$ 
{\mathcal I}(X) = \sum_{|d|\le X}^{\flat} L(\tfrac 12,\chi_d) \prod_{\ell =1}^{\infty} A_\ell(\chi_d), 
$$ 
where the sum is over fundamental discriminants $d$, and $A_\ell(\chi_d) = 
\sum_{n\le X^{\vartheta/a_\ell}} d_{k/a_\ell}(n) \chi_d(n)/\sqrt{n}$.  Arguing as in \cite{RuS2} 
we may evaluate ${\mathcal I}(X)$ for suitably small $\vartheta$, and then 
obtain a lower bound for $\sum_{|d|\le X}^{\flat} |L(\tfrac 12, \chi_d)|^k$.

\end{document}